\documentclass[12pt, reqno]{amsart}
\usepackage{amssymb,amsthm,amsfonts,amsmath, amscd}

\usepackage{hyperref}
\usepackage{mathrsfs}

\newcommand\A{\mathbb A}
\newcommand\Y{\mathbb Y}
\newcommand\Z{\mathbb Z}
\newcommand\C{\mathbb C}
\newcommand\R{\mathbb R}

\newcommand\al{\alpha}
\newcommand\be{\beta}
\newcommand\ga{\gamma}
\newcommand\Ga{\Gamma}
\newcommand\de{\delta}

\newcommand\La{\Lambda}
\newcommand\la{\lambda}

\newcommand\epsi{\varepsilon}
\newcommand\om{\omega}
\newcommand\Om{\Omega}

\newcommand\MM{\mathfrak M}

\newcommand\LL{\mathfrak L}

\newcommand\wt{\widetilde}

\newcommand\Conf{\operatorname{Conf}}

\newcommand\Fun{{\operatorname{Fun}}}

\newcommand\ME{{\mathsf M}}

\newcommand\D{\mathfrak D}

\newcommand\ext{{\operatorname{ext}}}

\newcommand\pd{\partial}
\newcommand\zz{{z,z'}}
\newcommand\zxi{{z,z',\xi}}

\newcommand\FS{F\!S}

\newcommand\ord{{\operatorname{ord}}}

\newtheorem{theorem}{Theorem}[section]

\newtheorem{corollary}[theorem]{Corollary}

\theoremstyle{definition}
\newtheorem{definition}[theorem]{Definition}
\newtheorem{remark}[theorem]{Remark}

\numberwithin{equation}{section}

\begin{document}

\title[Laguerre and Meixner symmetric functions]
{Laguerre and Meixner symmetric functions,\\ and infinite-dimensional diffusion
processes}

\author{Grigori Olshanski}
\address{Institute for Information Transmission Problems, Bolshoy
Karetny 19,  Moscow 127994, Russia; Independent University of Moscow, Russia}
\email{olsh2007@gmail.com}

\thanks{Supported by the RFBR grant 08-01-00110, the RFBR-CNRS grant 10-01-93114,
and the project SFB 701 of Bielefeld University.}

\begin{abstract}
The Laguerre symmetric functions introduced in the note are indexed by
arbitrary partitions and depend on two continuous parameters. The top degree
homogeneous component of every Laguerre symmetric function coincides with the
Schur function with the same index. Thus, the Laguerre symmetric functions form
a two-parameter family of inhomogeneous bases in the algebra of symmetric
functions. These new symmetric functions are obtained from the $N$-variate
symmetric polynomials of the same name by a procedure of analytic continuation.
The Laguerre symmetric functions are eigenvectors of a second order
differential operator, which depends on the same two parameters and serves as
the infinitesimal generator of an infinite-dimensional diffusion process
$X(t)$. The process $X(t)$ admits approximation by some jump processes related
to one more new family of symmetric functions, the Meixner symmetric functions.

In equilibrium, the process $X(t)$ can be interpreted as a time-dependent point
process on the punctured real line $\R\setminus\{0\}$, and the point
configurations may be interpreted as doubly infinite collections of particles
of two opposite charges and log-gas-type interaction. The dynamical correlation
functions of the equilibrium process have determinantal form: they are given by
minors of the so-called extended Whittaker kernel, introduced earlier in a
paper by Borodin and the author.

\end{abstract}

\date{}

\maketitle

\tableofcontents

\section{Introduction}

\subsection{Preface}

The present note is a research announcement; the detailed exposition will
appear elsewhere. The goal of the work is twofold: (1) to introduce  new bases
$\{\LL_\nu\}$ and $\{\MM_\nu\}$ in the algebra $\La$ of symmetric functions,
and (2) to construct a diffusion process $X(t)$ in an infinite-dimensional cone
$\wt\Om$. The two subjects are interrelated: The algebra $\La$ serves as the
algebra of ``polynomial observables'' on the cone $\wt\Om$, and the basis
elements $\LL_\nu\in\La$ are the eigenfunctions of the infinitesimal generator
of the process $X(t)$. As for the basis $\{\MM_\nu\}$, its elements are the
eigenfunctions of the infinitesimal generator of an auxiliary jump process.

The basis elements $\LL_\nu\in\La$ are called the {\it Laguerre symmetric
functions\/}. They are indexed by arbitrary partitions $\nu$ and depend on two
continuous parameters. The diffusion process $X(t)$ depends on the same
parameters. It possesses an invariant symmetrizing measure $W$, which also
serves as the orthogonality measure for the Laguerre symmetric functions.

\subsection{Finite-dimensional counterparts}

All the basic objects --- the Laguerre symmetric functions, the cone $\wt\Om$,
the probability distribution $W$ on $\wt\Om$, and the diffusion $X(t)$ on
$\wt\Om$
--- have finite-dimensional counterparts; I will describe them briefly.

\medskip

$\bullet$ In dimension 1, these are the classical Laguerre orthogonal
polynomials on the half-line $x>0$ with the weight measure $x^{b-1}e^{-x}dx$
(here $b>0$ is a parameter), and the diffusion is generated by the associated
ordinary differential operator
$$
x\frac{d^2}{dx^2}+(b-x)\frac{d}{dx};
$$
the Laguerre polynomials are its eigenfunctions.

$\bullet$ In dimension $N=2,3,\dots$, we deal with the algebra $\La_N$ of
symmetric polynomials in $N$ variables $x_1,\dots,x_N$. Such polynomials are
viewed as functions on the $N$-dimensional cone
$$
\wt\Om_N=\{\underline{x}:=(x_1,\dots,x_N): x_1\ge\dots\ge x_N\ge0\}\subset\R^N.
$$
A relevant basis in $\La_N$ is formed by the $N$-variate symmetric Laguerre
polynomials \cite{Ma87}, \cite{La91c}, which are orthogonal with respect to the
measure on $\wt\Om_N$ with the density
$$
(x_1\dots x_N)^{b-1}e^{-x_1-\dots-x_N}\cdot\prod_{1\le i<j\le N}(x_i-x_j)^2.
$$
Assuming $\underline{x}$ to be in the interior of the cone, one can interpret
$\underline{x}$ as a collection of $N$ indistinguishable particles on the
half-line $\R_{>0}$; then the above measure determines an ensemble of random
particle configurations, called the $N$-particle {\it Laguerre ensemble\/}.
Again, there exists an associated diffusion process $X_N(t)$ with state space
$\wt\Om_N$. In the interior of the cone, $X_N(t)$ may be interpreted as a
random evolution of $N$ interacting particles with pairwise repulsion. One may
call $X_N(t)$ the $N$-particle {\it dynamical Laguerre ensemble\/}.

$\bullet$ Note also that there exists a lattice version of the Laguerre
ensemble, related to natural discrete analogs of the Laguerre polynomials ---
the Meixner polynomials. There exists also an associated Markov dynamics, which
is a Markov jump process.  In the simplest case $N=1$ this is a birth-death
process with linear jump rates.

\subsection{Analytic continuation}

In the literature, there exist many models of such a kind, continuous and
discrete, associated with various systems of orthogonal polynomials. A general
recipe for building infinite-dimensional analogs of such models, often applied
in Random Matrix Theory, is to use a large-$N$ limit transition (see, e.g., the
survey paper \cite{KT10c}). However, in the present work, a different approach
is applied. In short, its main idea can be formulated as {\it extrapolation
into complex domain via analytic continuation with respect to two parameters,
the number of particles $N=1,2,3,\dots$ and the additional parameter $b>0$\/}.
Surprisingly enough, although parameters $N$ and $b$ are of a very different
origin, they can be treated on equal grounds.

In comparison with the existing approaches to random-matrix-type dynamical
models in infinite dimension (see \cite{Sp87}, \cite{Os09}, \cite{KT09},
\cite{KT10a}, \cite{KT10b}), our approach is to great extent more algebraic.

\subsection{Point processes}

In the equilibrium state corresponding to the stationary distribution $W$, the
process $X(t)$ can be interpreted as a time-dependent point process
$X^{\operatorname{stat}}(t)$. This interpretation relies on the fact that the
invariant measure $W$ is supported by a dense Borel subset
$\wt\Om'\subset\wt\Om$ admitting a natural realization as a space of infinite
particle configuration on the punctured real line $\R^*:=\R\setminus\{0\}$. It
turns out that the time-dependent point process $X^{\operatorname{stat}}(t)$ is
{\it determinantal\/}: its dynamical correlation functions are given by minors
of a kernel $K(s,u; t,v)$ on $(\R\times\R^*)\times(\R\times\R^*)$ --- the
so-called {\it extended Whittaker kernel\/}, which initially appeared in
\cite{BO06a}.

Note that the $N$-particle dynamical ensembles $X_N(t)$ live on the half-line
$\R_{>0}$ while the particle configurations of the process
$X^{\operatorname{stat}}(t)$ live on the punctured real line $\R^*$, which is
the union of {\it two\/} half-lines. Because of this duplication effect, the
claim that $X(t)$ cannot be related to the processes $X_N(t)$ through a
large-$N$ limit transition becomes intuitively evident.

\subsection{Lattice approximation}

Although $X(t)$ does not arise from a large-$N$ limit, it admits a {\it lattice
approximation\/}. Namely, $X(t)$ can be obtained as a scaling limit of some
jump processes depending on an additional parameter $\xi\in(0,1)$, as
$\xi\uparrow1$. These jump processes were studied in detail in \cite{BO06a}.
Their state space is the set of all partitions. This countable set also can be
realized as a set of particle configurations on the lattice $\Z+\tfrac12$ of
half-integers; here the number of particles is finite but not restricted. The
second basis $\{\MM_\nu\}\subset\La$ mentioned above just arises in connection
with these jump processes. The elements $\MM_\nu$ are called the {\it Meixner
symmetric functions\/}; they depend on three parameters: to the two parameters
of the Laguerre symmetric functions one adds the third parameter $\xi\in(0,1)$.
In the limit as $\xi\uparrow1$, one has $\MM_\nu\to\LL_\nu$, similarly to the
approximation of the classical Laguerre polynomials by the Meixner polynomials.

\subsection{Concluding remarks}

The results announced in the present note continue those of \cite{BO00},
\cite{BO06a}, \cite{BO06b}, \cite{BO09}, and all these works have a connection
with the asymptotic representation theory of the symmetric groups. It is
interesting to compare the results of these papers and the present note with
those of the papers \cite{BO05a}, \cite{BO05b}, \cite{BO10}, which are related
to representations of the unitary groups. Although the construction of a Markov
dynamics in \cite{BO10} relies on a different approach, \cite[Appendix]{BO10}
also exploits the idea of analytic continuation. Note that in the context of
the unitary groups, the Laguerre polynomials are (to some extent) replaced by
the Jacobi polynomials.

Finally, I would like to note that the representation theory of reductive
groups and Lie algebras also affords examples in which finite-dimensional
objects arise as a degeneration of infinite-dimensional ones, and, conversely,
infinite-dimensional objects may be reconstructed from finite-dimensional ones
through analytic continuation in parameters. For instance, the principal series
representations or highest weight modules may be viewed as analytic
continuation of the irreducible finite-dimensional representations.

\section{The one-particle case}\label{s2}

There is a well-known relationship between systems $\{\phi_n; n=0,1,2,\dots\}$
of orthogonal polynomials of hypergeometric type on $\R$ and some
one-dimensional Markov processes $X(t)$ (see, e.g., \cite{Sch00}). Namely, the
state space of $X(t)$ is a closed subset $I\subset\R$ --- the support of the
orthogonality measure of $\{\phi_n\}$, and the infinitesimal generator of
$X(t)$ is a second order differential (or difference) operator $D$, such that
$$
D \phi_n=-\mu_n \phi_n, \qquad \mu_n\ge0.
$$
The orthogonality measure $w$ for the polynomials $\phi_n$ serves as an
invariant and symmetrizing measure of $X(t)$. The transition function of $X(t)$
has the form
\begin{equation}\label{eq4}
P(t;x,dy)=\left(\sum_{n=0}^\infty e^{-\mu_nt}\frac{\phi_n(x)\phi_n(y)}{\int_I
\phi^2_n(u)w(du)}\right)w(dy)
\end{equation}

The simplest examples are provided by the classical Jacobi, Laguerre, and
Hermite polynomials:

$\bullet$ {\it Jacobi polynomials\/}: $I$ is the closed interval $[-1,1]$, $w$
has density $(1-x)^{a-1}(1+x)^{b-1}$ with parameters $a,b>0$, and
$$
D=(1-x^2)\frac{d^2}{dx^2}+[b-a-(a+b)x]\frac{d}{dx}.
$$

$\bullet$ {\it Laguerre polynomials\/}: $I$ is the closed half-line
$[0,+\infty)$, $w$ has density $x^{b-1}e^{-x}$ with parameter $b>0$, and
$$
D=x\frac{d^2}{dx^2}+(b-x)\frac{d}{dx}.
$$

$\bullet$ {\it Hermite polynomials\/}: $I$ is the whole real line, $w$ has
density $e^{-x^2/2}$, and
$$
D=\frac{d^2}{dx^2}-x\frac{d}{dx}.
$$

\medskip

In the Hermite case, $X(t)$ is the Ornstein-Uhlenbeck process. In the Laguerre
case, $X(t)$ is closely related to a squared Bessel process (see, e.g.
\cite{Eie83}).

\section{The $N$-particle case}

Here we recall a well-known construction providing a multidimensional
generalization of the above picture.

Fix $N=1,2,3,\dots$\,. Instead of univariate polynomials we will deal with {\it
symmetric\/} polynomials in $N$ variables $x_1,\dots,x_N$. Denote by $\La_N$
the algebra of such polynomials (the base field is $\R$ or $\C$ depending on
convenience). The interval $I$ is replaced by the subset
$$
I^N_\ord:=\{(x_1,\dots,x_N)\subset I^N: x_1\ge\dots\ge x_N\},
$$
and we regard $\La_N$ as an algebra of functions on $I^N_\ord$. Let
$\nu=(\nu_1,\dots,\nu_N)$ range over the set of partitions of length
$\ell(\nu)\le N$. We set
\begin{equation}\label{eq1}
\phi_{\nu\mid
N}(x_1,\dots,x_N)=\frac{\det[\phi_{\nu_i+N-i}(x_j)]}{V_N(x_1,\dots,x_N)},
\end{equation}
where the determinant is of order $N$ and $V_N$ is the Vandermonde,
$$
V_N=V_N(x_1,\dots,x_N)=\prod_{1\le i<j\le N}(x_i-x_j).
$$
The $\phi_{\nu\mid N}$'s are symmetric polynomials forming a basis in $\La_N$.
Moreover, it is readily verified that they are pairwise orthogonal with respect
to the measure
\begin{equation}\label{eq2}
w_N(dx_1\dots dx_N)=(V_N)^2\prod_{i=1}^Nw(dx_i)
\end{equation}
on $I^N_\ord$.

This construction seems to be well known; see e.g. Lassalle's papers
\cite{La91a}, \cite{La91b}, \cite{La91c}. Formula \eqref{eq1} is similar to the
classical expression for the Schur symmetric polynomials (the Schur polynomials
appear if one substitutes $\phi_n(x)=x^n$; they are not orthogonal polynomials
though).

Further, the analog of $D$ is the operator
\begin{equation}\label{eq3}
D_N:=V_N^{-1}\circ(D^{x_1}+\dots+ D^{x_N})\circ V_N+d_N1
\end{equation}
where $D^{x_i}$ stands for a copy of $D$ acting on variable $x_i$ and
$$
d_N=\mu_0+\dots+\mu_{N-1}.
$$

If $D$ is a differential operator,
$$
D=A(x)\frac{d^2}{dx^2}+B(x)\frac{d}{dx},
$$
then $D_N$ is a partial differential operator,
$$
D_N:=\sum_{i=1}^N\left(A(x_i)\frac{\pd^2}{\pd x_i^2}+\left[B(x_i)+\sum_{j:\,
j\ne i}\frac{2A(x_i)}{x_i-x_j}\right]\frac{\pd}{\pd x_i}\right).
$$
Due to the special choice of $d_N$, the constant term in $D_N$ vanishes.

Although the coefficients in front of the first derivatives have singularities
along the hyperplanes $x_i=x_j$, the operator $D_N$ is applicable to symmetric
polynomials and preserves the space $\La_N$. The polynomials $\phi_{\nu\mid N}$
are eigenfunctions of $D_N$:
$$
D_N \phi_{\nu\mid
N}=-\left(\sum_{i=1}^N(\mu_{\nu_i+N-i}-\mu_{N-i})\right)\phi_{\nu\mid N}.
$$

Finally, $D_N$ serves as the (pre-)generator of a Markov process $X_N(t)$ on
$I^N_\ord$ with invariant symmetrizing measure $w_N$.

The first example of such a process $X_N(t)$ has been investigated in
\cite{Dy62}; it corresponds to the system of Hermite polynomials. As shown in
that paper, $X_N(t)$ is obtained from a matrix-valued Ornstein-Uhlenbeck
process through the projection onto the matrix eigenvalues. One may say that
$X_N(t)$ is the radial part of this matrix-valued Markov process.

In the present work we focus on the Laguerre case.

\section{The Laguerre symmetric functions}

Let $e_1,e_2,\dots$ denote the elementary symmetric polynomials,
$$
e_1=\sum_i x_i, \quad e_2=\sum_{i<j}x_i x_j, \quad e_3=\sum_{i<j<k}x_i x_j x_k,
$$
and so on. Here it is tacitly assumed that the indices range over
$\{1,\dots,N\}$, where $N$ is the number of variables. As well known, the
algebra $\La_N$ of $N$-variate symmetric polynomials is isomorphic to the
algebra of ordinary polynomials in $e_1,\dots,e_N$.

Our first step is to make a change of variables: take $\{e_1,\dots,e_N\}$ as
new (formal) variables instead of natural coordinates $x_1,\dots,x_N$.

\begin{theorem}
The $N$-variate Laguerre operator $D_N:\La_N\to\La_N$ can be rewritten as the
following second order partial differential operator in variables
$e_1,\dots,e_N$:
$$
D_N=\sum_{m,n=1}^N A_{mn}\frac{\pd^2}{\pd e_m\pd e_n}+\sum_{n=1}^N B_n
\frac{\pd}{\pd e_n},
$$
where
$$
A_{mn}=\sum_{k=0}^{\min(m,n)-1}(m+n-1-2k)e_{m+n-1-k}e_k
$$
and
$$
B_n=-n e_n+(N-n+1)(N+b-n)e_{n-1}
$$
with the convention that $e_0=1$ and $e_k=0$ for $k>N$.
\end{theorem}

The next step is to replace $\La_N$ by the {\it algebra $\La$ of symmetric
functions\/}. For our purpose, it is convenient to define $\La$ as the algebra
of polynomials in countably many  formal commuting variables $e_1,e_2,\dots$,
which are assumed to be algebraically independent.

For $N=1,2,\dots$, let $J_N\subset\La$ denote the ideal generated by elements
$e_k$ with $k>N$. The quotient algebra $\La/J_N$ is naturally isomorphic to
$\La_N$, so we get a canonical algebra homomorphism $\pi_N:\La\to\La_N$, which
we call the $N$th {\it truncation map\/}.

\begin{definition}\label{dfnB}
Let $z$ and $z'$ be complex parameters. Consider the formal second order
differential operator in countably many variables $e_1,e_2,\dots$, obtained
from the $N$-variate Laguerre operator $D_N$ by removing the relations
$e_{N+1}=e_{N+2}=\dots=0$, dropping the restriction $m,n\le N$, and replacing
the factor $(N-n+1)(N+b-n)$ in the definition of coefficient $B_n$ by
$(z-n+1)(z'-n+1)$:
$$
\D=\sum_{m,n=1}^\infty A_{mn}\frac{\pd^2}{\pd e_m\pd e_n}+\sum_{n=1}^\infty B_n
\frac{\pd}{\pd e_n},
$$
where $A_{mn}$ is given by exactly the same formula as above,
$$
A_{mn}=\sum_{k=0}^{\min(m,n)-1}(m+n-1-2k)e_{m+n-1-k}e_k,
$$
and
$$
B_n=-n e_n+(z-n+1)(z'-n+1)e_{n-1}.
$$
Observe that $\D$ is correctly defined as an operator $\La\to\La$; we call it
the {\it Laguerre operator in\/ $\La$.}
\end{definition}

If $z=N$ and $z'=N+b-1$, then $\D$ preserves the ideal $J_N\subset\La$ and
hence factorizes to an operator in the quotient $\La/J_N=\La_N$; the resulting
operator coincides with $D_N$. Note that this property, combined with the
polynomial dependence on parameters $z,z'$, determines the operator uniquely.
In this sense, $\D$ may be viewed as the result of {\it analytic
continuation\/} (or extrapolation) of the $N$-variate Laguerre operators $D_N$
with respect to parameters $N$ and $b$.

Let $\{L_n\}$ denote the system of {\it monic\/} Laguerre polynomials with
parameter $b>0$ (see, e.g., \cite{KS96}). Recall that, in our notation, the
weight function is $x^{b-1}e^{-x}$. Next, let $\{L_{\nu\mid N,b}\}$ denote the
system of $N$-variate symmetric Laguerre polynomials defined in accordance with
the determinantal formula \eqref{eq1}. We are going to define elements of $\La$
that may be viewed as are analogs of the polynomials $L_{\nu\mid N,b}$. To do
this we apply the same principle of analytic continuation in $N$ and $b$ as we
have employed in the definition of $\D$.

\begin{theorem}\label{thmB}
For an arbitrary partition $\nu=(\nu_1,\nu_2,\dots)$, there exists a unique
element $\LL_\nu\in\La[z,z']:=\La\otimes\C[z,z']$, such that for any natural
number $N\ge\ell(\nu)$ and any $b>0$,
$$
\pi_N\left(\LL_\nu\big|_{z=N,\,z'=N+b-1}\right)=L_{\nu\mid N,b}.
$$
\end{theorem}

Here $\pi_N:\La\to\La_N$ is the $N$th truncation map introduced above.

\begin{definition}
We call the elements $\LL_\nu$ the {\it Laguerre symmetric functions\/}.
\end{definition}

Recall the definition of the {\it Schur symmetric functions\/}: these are
elements of $\La$ indexed by arbitrary partitions $\nu$ and expressed through
the generators $e_n$ by the following formula (the {\it N\"agelsbach--Kostka
formula\/}, see \cite{Ma95})
$$
S_\nu=\det[e_{\nu\,'-i+j}];
$$
here $\nu\,'$ stands for the partition given by transposing the Young diagram
corresponding to $\nu$, and the order of the determinant is an arbitrary
integer $\ge\ell(\nu\,')$. As well known, the Schur functions form a
distinguished homogeneous basis in $\La$, and
$$
\deg S_\nu=|\nu|:=\sum \nu_i.
$$

Let us explain the notation used in the next theorem. We identify partitions
and the corresponding Young diagrams. Given a couple $\mu\subseteq\nu$ of Young
diagrams, we denote by $\dim\nu/\mu$ the number of standard Young tableaux of
the skew shape $\nu/\mu$. The symbol $\Box\in\nu/\mu$ denotes a box in
$\nu/\mu$, and $c(\Box)$ denotes its {\it content\/}, equal to the difference
$j-i$ of the column number $j$ and the row number $i$ of the box.

\begin{theorem}\label{thmC}
The expansion of the  Laguerre symmetric function $\LL_\nu$ in the basis of the
Schur symmetric functions has the form
$$
\LL_\nu=\sum_{\mu\subseteq\nu}C(\nu,\mu;z,z')S_\mu,
$$
where
$$
C(\nu,\mu;z,z')=(-1)^{|\nu|-|\mu|} \frac{\dim\nu/\mu}{(|\nu|-|\mu|)!}\,
\prod_{\Box\in\nu/\mu}(z+c(\Box))(z'+c(\Box))).
$$
\end{theorem}

Since $C(\nu,\nu;z,z')=1$, the top homogeneous component of $\LL_\nu$ is equal
to $S_\nu$:
$$
\LL_\nu=S_\nu+\textrm{lower degree terms}.
$$
It follows that the Laguerre symmetric functions with any fixed values of
parameters $z$ and $z'$ form a basis in $\La$.

For the empty diagram corresponding to the zero partition, $\nu=\varnothing$,
we have $\LL_\varnothing=S_\varnothing=1$. This is the only case when $\LL_\nu$
and $S_\nu$ coincide: for $\nu\ne\varnothing$, $\LL_\nu$ is an inhomogeneous
element, so that the basis $\{\LL_\nu\}$ of Laguerre symmetric functions is an
example of inhomogeneous basis. In this respect, it differs from other bases in
$\La$ that are commonly used in algebraic combinatorics.

A box $\square$ in a Young diagram $\nu$ is said to be a {\it corner box\/} if
the shape $\nu\setminus\square$ obtained by removing $\square$ from $\nu$ is
again a Young diagram. Let $\nu^-$ denote the set of all corner boxes in $\nu$.
For instance, if $\nu=(3,2,2)$ then $\nu^-$ comprises two corner boxes,
$\square=(1,3)$ and $\square=(3,2)$.

\begin{theorem}\label{thmD}
The action of $\D$ on the Schur functions is given by
$$
\D S_\nu=-|\nu|S_\nu+\sum_{\Box\in\nu^-}
(z+c(\square))(z'+c(\Box))S_{\nu\setminus\Box}
$$
\end{theorem}

\begin{theorem}\label{thmE}
The Laguerre symmetric functions are eigenvectors of the Laguerre operator $\D$
with the same values of parameters $(z,z')$. More precisely,
$$
\D\LL_\nu=-|\nu|\LL_\nu.
$$
\end{theorem}

\section{Formal orthogonality}

\begin{definition}\label{dfnA}
For any fixed $(z,z')\in\C^2$, introduce the {\it formal moment functional\/}
$\psi:\La\to\C$ by setting
$$
\psi(1)=1, \qquad \psi(\D f)=0 \quad\textrm{for any $f\in\La$},
$$
where the Laguerre operator $\D$ is taken with the same values of the
parameters as $\psi$.
\end{definition}

The definition is correct by virtue of the last theorem. Indeed, it implies
that that the range of $\D$ is the span of the Laguerre functions $\LL_\nu$
with $\nu\ne\varnothing$, while the vector $1=\LL_\varnothing\in\La$ is
transversal to this span. Note that $\psi$ depends polynomially on the
parameters $z,z'$.

\begin{theorem}\label{thmF}
For any Young diagram $\nu$,
$$
\psi(S_\nu)=\prod_{\Box\in\nu}(z+c(\Box))(z'+c(\Box))\cdot
\left(\frac{\dim\nu}{|\nu|!}\right)^2
$$
\end{theorem}

This formula provides an alternative (equivalent) way of introducing the moment
functional.

\begin{theorem}\label{thmA}
For any two Young diagrams $\nu$ and $\mu$,
\begin{equation}\label{eq8}
\psi\left(\LL_\nu\LL_\mu\right)=\de_{\nu\mu}\prod_{\Box\in\nu}(z+c(\Box))(z'+c(\Box)),
\end{equation}
where $\de_{\nu\mu}$ is Kronecker's delta.
\end{theorem}

This result shows that the Laguerre symmetric functions are pairwise orthogonal
with respect to the inner product in the space $\La$ defined by
\begin{equation}\label{eq9}
(f,g)=(f,g)_{z,z'}:=\psi(fg), \quad f,g\in\La.
\end{equation}
Obviously, the inner product is nondegenerate if (and only if) $z$ and $z'$ are
not integers, for then the product in the right-hand side of \eqref{eq8} never
vanishes.

\begin{remark}\label{rmk2}
Assume $z$ and $z'$ are not integers. Then the Laguerre symmetric function
$\LL_\nu$ is characterized by the following two properties:

(1) $\LL_\nu$ differs from the Schur symmetric function $S_\nu$ by lower degree
terms;

(2) $\LL_\nu$ is orthogonal, with respect to inner product \eqref{eq9}, to all
elements of $\La$ of lower degree (that is, of degree strictly less than
$\deg\LL_\nu=|\nu|$).

Alternatively, without any assumption on the parameters, $\LL_\nu$ is
characterized by (1) together with the following property replacing (2):

(2\,$'$) $\LL_\nu$ is an eigenfunction of the Laguerre operator $\D$.
\end{remark}

\section{The orthogonality measure for the Laguerre symmetric functions}

\begin{definition}
By the {\it Thoma cone\/} we mean the subset
$\wt\Om\subset\R_+^\infty\times\R_+^\infty\times\R_+$ consisting of triples
$\om=(\al,\be,\de)$, where
$$
\al=(\al_1\ge\al_2\ge\dots\ge0), \quad \be=(\be_1\ge\be_2\ge\dots\ge0), \quad
\de\ge0
$$
and
$$
\sum_{i=1}^\infty(\al_i+\be_i)\le\de.
$$
The {\it Thoma simplex\/} $\Om$ is the subset of $\wt\Om$ determined by the
additional condition $\de=1$.
\end{definition}

Both $\wt\Om$ and $\Om$ are closed subsets of the product space
$\R^{2\infty+1}:=\R^\infty\times\R^\infty\times\R$ equipped with the product
topology; $\Om$ is compact and $\wt\Om$ is locally compact; $\wt\Om$ is
precisely the cone with base $\Om$.

Let $\Fun(\wt\Om)$ denote the space of continuous functions on the Thoma cone;
this is an algebra under pointwise multiplication. We embed $\La$ into
$\Fun(\wt\Om)$ by setting
\begin{equation}\label{eq7}
1+\sum_{n=1}^\infty e_nt^n \, \to\, e^{\ga t}\prod_{i=1}^\infty\frac{1+\al_i
t}{1-\be_i t},
\end{equation}
where $t$ is an auxiliary formal variable and
$$
\ga:=\de-\sum_{i=1}^\infty(\al_i+\be_i).
$$
More precisely, the correspondence, which is written in terms of the generating
series for $\{e_n\}$, turns each $e_n$ into a function $e_n(\om)$ on the Thoma
cone. This function is continuous. The correspondence $e_n\mapsto
e_n(\,\cdot\,)$ is extended to the whole algebra $\La$ by multiplicativity. In
this way we get an algebra morphism $\La\to \Fun(\wt\Om)$, which is an
embedding. For any $f\in\La$ we denote the corresponding continuous function on
$\wt\Om$ by $f(\om)$.

Equivalently, in terms of another system of generators of the algebra $\La$,
the Newton power sums $p_n$, the embedding $\La\to \Fun(\wt\Om)$ can be defined
by setting
$$
p_1(\om)=\de, \qquad p_n(\om)=\sum_{i=1}^\infty
\al_i^n-\sum_{i=1}^\infty(-\be_i)^n \quad \textrm{for $n\ge2$}.
$$

In the algebra $\La$, there is a distinguished involutive automorphism, which
is defined on the generators $p_n$ as $p_n\mapsto (-1)^{n-1}p_n$. Under this
automorphism, $S_\nu\mapsto S_{\nu\,'}$. The above formula shows that in the
realization $\La\subset \Fun(\wt\Om)$, this automorphism amounts to
transposition $\al\leftrightarrow\be$.

\begin{definition}
Let us say that $(z,z')\in\C^2$ is {\it admissible\/} if both $z$ and $z'$ are
nonzero and
$$
\prod_{\Box\in\nu}(z+c(\Box))(z'+c(\Box))\ge0
$$
for any Young diagram $\nu$.
\end{definition}

The family of admissible couples $(z,z')$ splits into the union of the three
subfamilies:

$\bullet$ The {\it principal series\/}: both $z$ and $z'$ are nonreal, $z'=\bar
z$.

$\bullet$ The {\it complementary series\/}: both $z$ and $z'$ are real and are
contained inside an open interval $(m,m+1)$ with $m\in\Z$.

$\bullet$ The {\it degenerate series\/}: $(z,z')=\pm(N,N+b-1)$ or
$(z',z)=\pm(N, N+b-1)$, where $N=1,2,\dots$ and $b>0$.

\begin{theorem}\label{thmG}
Assume $(z,z')$ is admissible. There exists a probability measure $W=W_{z,z'}$
on the Thoma cone $\wt\Om$ such that all functions from $\La$ are
$W$-integrable and the formal moment functional $\psi$ with parameters $(z,z')$
coincides with expectation with respect to $W$:
$$
\psi(f)=\int_{\wt\Om}f(\om)W(d\om) \qquad\text{\rm for all $f\in\La$}.
$$
These properties determine $W$ uniquely. Moreover, the functions from $\La$ are
dense in the Hilbert space $L^2(\wt\Om, W)$.
\end{theorem}

By virtue of Theorem \ref{thmA} this implies that the measure $W$ with
admissible parameters $(z,z')$ is the orthogonality measure for the Laguerre
symmetric functions.

In the case when $z=N=1,2,\dots$ and $z'=N+b-1$ with $b>0$ we recover the
Laguerre measures $w_N$ on the $N$-dimensional cone
$$
\wt\Om_N:=\{(x_1\ge\dots x_N\ge0)\}.
$$
Here we embed $\wt\Om_N$ into the Thoma cone by setting $\al_i=x_i$ for
$i=1,\dots,N$ and $\de=\sum x_i$ (so that all remaining $\al$- and
$\be$-coordinates are equal to zero).

An immediate consequence of the theorem is

\begin{corollary}
Assume that the couple $(z,z')$ is admissible and nondegenerate, i.e. belongs
to the principal or complementary series. Then the Laguerre symmetric
functions, viewed as functions on the Thoma cone, form an orthogonal basis in
the Hilbert space $L^2(\wt\Om,W_{z,z'})$.
\end{corollary}

\section{Properties of the orthogonality measure}

Here are some properties of the measures $W_{z,z'}$ with admissible parameters
$(z,z')$:

$\bullet$ $W_{z,z'}$ are obtained from the so-called z-measures on the Thoma
simplex by a simple integral transform along the rays of the Thoma cone. See
\cite[Section 5]{BO00}.

$\bullet$ $W_{z,z'}$ does not change under transposition $z\leftrightarrow z'$.
Within this symmetry relation, the measures corresponding to different couples
of parameters are pairwise disjoint: this follows from \cite{KOV04}.

$\bullet$ The involutive map $\al\leftrightarrow\be$ of the Thoma cone
transforms $W_{z,z'}$ to $W_{-z,-z'}$.

$\bullet$ $W_{z,z'}$ is supported by the subset of the Thoma cone formed by
those triples $\om=(\al,\be,\de)$ for which $\de$ exactly equals
$\sum(\al_i+\be_i)$, i.e. $\ga=0$. This subset is Borel and everywhere dense in
$\wt\Om$. (Note that $\ga$ is not a continuous function in $\om$, it is only
lower semicontinuous.)

\begin{theorem}
If $(z,z')$ is in the principal or complementary series, then the topological
support of $W_{z,z'}$ is the whole space $\wt\Om$.
\end{theorem}

(The topological support of a measure is the smallest closed subset supporting
the measure.)

The theorem implies that a nonzero continuous function on $\wt\Om$, which is
square integrable with respect to $W_{z,z'}$, is a nonzero vector in the
Hilbert space $L^2(\wt\Om,W_{z,z'})$. In particular, the natural map assigning
to a bounded continuous function a vector in this Hilbert space is an
embedding. This assertion fails in the case when $(z,z')$ belongs to the
degenerate series.

\section{The Laguerre diffusion process on the Thoma cone}

{\it In this section $(z,z')$ is a fixed couple of parameters from the
principal or complementary series\/}.

Recall that we may regard $\La$ as a dense subspace in the Hilbert space
$H:=L^2(\wt\Om,W_{z,z'})$ and then $\{\LL_\nu\}$ becomes an orthogonal basis in
$H$. The Laguerre operator $\D:\La\to\La$ is diagonalized in this basis,
$$
\D\LL_\nu=-|\nu|\LL_\nu,
$$
so that the eigenvalues of $\D$ are $0,-1,-2,\dots$, where $0$ has multiplicity
$1$ and corresponds to the basis vector $\LL_\varnothing=1$. It follows that
$\D$ generates a strongly continuous semigroup $P(t)$ of contractive
selfadjoint operators in $H$,
$$
P(t)\LL_\nu=e^{-|\nu|t}\LL_\nu, \qquad t\ge0.
$$

\begin{theorem}
The semigroup $P(t)$ is a conservative Markov $L^2$-semigroup
\end{theorem}

By definition, this means that $P(t)1=1$ and $P(t)$ preserves the cone of
nonnegative functions in $H$. The first claim is obvious, the second one is
nontrivial; its proof relies on the approximation by some jump Markov
processes, see Section \ref{s14}.

The next claim says that $P(t)$ is actually a Feller semigroup (in one of the
versions of this property). Let $C(\wt\Om)\subset\Fun(\wt\Om)$ be the space of
bounded continuous functions on $\wt\Om$ with the supremum norm and
$C_0(\wt\Om)\subset C(\wt\Om)$ be its closed subspace formed by the functions
vanishing at infinity. Both $C(\wt\Om)$ and $C_0(\wt\Om)$ are Banach spaces
contained in $H$, but $C_0(\wt\Om)$ is separable while $C(\wt\Om)$ is not.

\begin{theorem}
The semigroup $P(t)$ preserves $C_0(\wt\Om)$ and induces a strongly continuous
contractive semigroup in this Banach space.
\end{theorem}

One of the ingredients of the proof is separation of variables described in the
next section.

The theorem implies that $P(t)$ gives rise to a Markov process $X(t)$ on
$\wt\Om$ with c\`adl\`ag sample trajectories. Actually, more can be said:

\begin{theorem}
With probability one, the trajectories of $X(t)$ are continuous.
\end{theorem}

By the very construction of the Markov process $X(t)$, the probability measure
$W_{z,z'}$ is its invariant and symmetrizing measure.

\section{Separation of variables}

Extend the algebra $\La$ by allowing division by $e_1$ (that is, localize over
the multiplicative semigroup generated by $e_1$) and denote the resulting
algebra by $\La^\ext$:
$$
\La^\ext:=\C[e_1, e_1^{-1}; e_2,e_3,\dots]\supset\La.
$$
Since the coefficients of the Laguerre differential operator $\D:\La\to\La$
(Definition \ref{dfnB} are polynomials in variables $e_n$, $\D$ can be extended
to an operator $\D^\ext:\La^\ext\to\La^\ext$.

Set
$$
r:=e_1, \quad e^\circ_n:=e_ne_1^{-n}, \qquad n\ge2,
$$
and
$$
e^\circ_0=e^\circ_1=1.
$$
The algebra $\La^\ext$ can be identified with $\C[r, r^{-1};
e^\circ_2,e^\circ_3,\dots]$.

\begin{theorem}
Under this identification, the operator $\D^\ext:\La^\ext\to\La^\ext$ takes the
form
$$
\D^\ext=\left(r\frac{\pd^2}{\pd r^2}+(c-r)\frac{\pd}{\pd
r}\right)+\frac1{r}\D^\circ,
$$
where
$$
\qquad c:=zz'
$$
and
$$
\D^\circ=\sum_{m,n\ge2}A^\circ_{mn}\frac{\pd^2}{\pd e^\circ_m\pd
e^\circ_n}+\sum_{n\ge1}B^\circ_n\frac{\pd}{\pd e^\circ_n}
$$
with
$$
A^\circ_{mn}=-mne^\circ_m
e^\circ_n+\sum_{k=0}^{\min(m,n)-1}(m+n-1-2k)e^\circ_{m+n-1-k}e^\circ_k
$$
and
$$
B^\circ_n=-n(n-1+c)e^\circ_n+(z-n+1)(z'-n+1)e^\circ_{n-1}.
$$
\end{theorem}

This result shows (at least on algebraic level) that the process $X(t)$ is the
skew product of the one-dimensional Laguerre diffusion with parameter $c=zz'>0$
and a Markov process on the Thoma simplex $\Om$, which is generated by the
operator $\D^\circ$. This can be compared to the splitting of the
multidimensional Brownian motion into the skew product of a one-dimensional
diffusion (a Bessel process) and the spherical Brownian motion, see, e.g.,
\cite{IM65}.

In more detail: The algebra $\La^\ext$ is realized, in a natural way, as an
algebra of functions on $\wt\Om\setminus\{0\}$. In this realization, elements
of the subalgebra $\La^\circ$ turn into homogeneous functions of degree $0$
with respect to homotheties of the cone $\wt\Om$. Since the Thoma simplex $\Om$
is a base of the Thoma cone, we may regard $\La^\circ$ as an algebra of
functions on $\Om$. Note that $\La^\circ$ is dense in the Banach space $C(\Om)$
of continuous functions on $\Om$. Thus, $\D^\circ$ becomes a densely defined
operator in $C(\Om)$. As shown in \cite{BO09}, the closure of $\D^\circ$
generates a diffusion process $X^\circ(t)$ in the Thoma simplex. \footnote{In
\cite{BO09}, the operator $\D^\circ$ is written down in a different coordinate
system.} Continuing the analogy with the multidimensional Brownian motion one
can say that in this picture, $X^\circ(t)$ is a counterpart of the spherical
Brownian motion and the variables $e^\circ_n$ play the role of spherical
coordinates.

\section{Correlation functions}

Set $\R^*=\R\setminus\{0\}$ and let $\Conf(\R^*)$ be the space of locally
finite point configurations on $\R^*$. Following \cite{BO00}, we define a
projection $\wt\Om\to\Conf(\R^*)$ by setting
$$
\om=(\al,\be,\de)\mapsto\{\al_i: \al_i\ne0\}\sqcup\{-\be_i:\be_i\ne0\}.
$$
Two points of $\wt\Om$ are mapped to one and the same configuration if and only
if they differ only by the value of the coordinate $\de$. Consequently, the
restriction of the projection on the subset
$$
\wt\Om':=\{\om: \de=\sum(\al_i+\be_i)\}
$$
is injective.

The {\it equilibrium version\/} $X^{\operatorname{stat}}(t)$ of the process
$X(t)$ is obtained by taking the stationary distribution $W$ as the initial
one. The process $X^{\operatorname{stat}}(t)$ is stationary in time; moreover,
since $W$ is a symmetrizing measure, one may extend the time parameter $t$ from
the half-line $[0,+\infty)$ to the whole real line $\R$. We know that the
stationary distribution $W$ is concentrated on the subset $\wt\Om'$. It follows
that the finite-dimensional distributions of the equilibrium process may be
viewed as probability measures on the spaces $(\wt\Om')^k$, $k=1,2,\dots$\,.
Next, applying the above projection, we may interpret every $k$-dimensional
distribution ($k=1,2,\dots$) as a probability measure on
$$
(\Conf(\R^*))^k=\Conf(\underbrace{\R^*\sqcup\dots\sqcup\R^*}_{\textrm{$k$
times}}),
$$
which is again a space of configurations. Such measures can be described in
terms of the correlation functions. In other words, these are the {\it
dynamical\/} or {\it space-time\/} correlation functions of a time-dependent
point process.

\begin{theorem}
The space-time correlation functions of the equilibrium process
$X^{\operatorname{stat}}(t)$ are determinantal. That is, they are given by
minors of a kernel $K(s,x;t,y)$, where $s$ and $t$ are time variables and
$x,y\in\R^*$.
\end{theorem}

The correlation kernel $K(s,x;t,y)$ is called the {\it extended Whittaker
kernel\/}. It was derived in \cite[Theorem B]{BO06a} as a scaling limit of the
correlation kernels of some equilibrium Markov jump processes on partitions.
Explicit expressions for $K(s,x;t,y)$ are contained in that paper.

Note, however, that the paper \cite{BO06a} left open the question whether the
kernel $K(s,x;t,y)$ determines a {\it Markov\/} process (a subtlety here is
that, in principle, it may happen that the Markov property is destroyed in a
limit transition). The results of the present section settle this question in
the positive.

\section{The Meixner symmetric functions}\label{s11}

Let $\Z_+\subset\Z$ denote the set of nonnegative integers. Fix two parameters
$b$ and $\xi$, where  $b>0$ as before and $0<\xi<1$.

The {\it classical Meixner polynomials\/} $M_n(x)$ are defined as the
orthogonal polynomials corresponding to the following discrete probability
measure supported by $\Z_+$:
$$
w^\ME=(1-\xi)^b\sum_{x\in\Z_+}\frac{(b)_x}{x!}\xi^x\de_x,
$$
where $\de_x$ denotes the Dirac measure at $x$ (see, e.g., \cite{KS96}). The
measure $w^\ME$ is is known under the name of the {\it negative binomial
distribution\/}. We use the standardization in which the $M_n$'s are monic
polynomials: $M_n(x)=x^n+\dots$.

Consider the following second order difference operator on $\Z_+$:
\begin{equation*}
\begin{aligned}
D^\ME f(x)& =\frac{\xi(b+x)}{1-\xi}f(x+1)+\frac x{1-\xi}f(x-1)\\
&-\frac{\xi(b+x)+x}{1-\xi}f(x)
\end{aligned}
\end{equation*}
(the factor $1-\xi$ in the denominator is introduced to simplify some formulas
below). $D^\ME$ is formally symmetric with respect to the weight function and
annihilates the constants. The Meixner polynomials are eigenfunctions of this
operator,
\begin{equation*}
D^\ME M_n=-n M_n.
\end{equation*}
Moreover, they can be characterized as the only polynomial eigenfunctions of
$D^\ME$.

The {\it $N$-variate symmetric Meixner polynomials\/}
$$
M_{\nu}(x_1,\dots,x_N)=M_{\nu\mid N, b,\xi}(x_1,\dots,x_N)
$$
are introduced following the recipe \eqref{eq1}.

Set
$$
\Z^N_{+,\ord}=\{(x_1,\dots,x_N)\in\Z^N_+: x_1>\dots>x_N\}
$$
and regard polynomials from $\La_N$ as functions on $\Z^N_{+\ord}$. Then the
Meixner polynomials $M_\nu$ become orthogonal polynomials with respect to the
atomic measure $w^\ME_N$ on $\Z^N_{+,\ord}$ defined according to \eqref{eq2}.

The polynomials $M_\nu$ are eigenfunctions of an operator
$D^\ME_N:\La_N\to\La_N$, which is defined according to \eqref{eq3}:
$$
D^\ME_N M_\nu=-|\nu|M_\nu.
$$
This operator can be realized as a difference operator on $\Z^N_{+\ord}$ acting
on a function $f$ by
$$
\begin{aligned}
D^\ME_Nf(x)&=\sum_{i=1}^N A_i(x)f(x+\epsi_i)+\sum_{i=1}^N B_i(x)f(x-\epsi_i) -
C(x)f(x)\\
&=\sum_{i=1}^N A_i(x)(f(x+\epsi_i)-f(x))+\sum_{i=1}^N
B_i(x)(f(x-\epsi_i)-f(x)).
\end{aligned}
$$
Here $x=(x_1,\dots,x_N)$ ranges over $\Z^N_{+,\ord}$,
$\{\epsi_1,\dots,\epsi_N\}$ is the canonical basis in $\R^N$, and the
coefficients are given by
$$
\begin{aligned}
A_i(x)&=\frac{\xi(b+x_i)}{1-\xi}\prod_{j:\,j\ne i}\frac{x_i-x_j+1}{x_i-x_j},\\
B_i(x)&=\frac{x_i}{1-\xi}\prod_{j:\,j\ne i}\frac{x_i-x_j-1}{x_i-x_j}\\
C(x)&=\frac{\xi b N+(1+\xi)\sum_{i=1}^Nx_i}{1-\xi}-\frac{N(N-1)}2.
\end{aligned}
$$

We will need a modified version of the truncation map $\pi_N:\La\to\La_N$; this
is an algebra morphism $\pi'_N:\La\to\La_N$, which is defined on the generators
$p_k$ (Newton power sums) by
$$\pi'_N(p_k)(x_1,\dots,x_N)=\sum_{i=1}^N[(x_i-N+\tfrac12)^k-(-i+\tfrac12)^k].
$$
Since the right-hand side is a symmetric polynomial, the definition makes
sense. It can be better understood in terms of the realization
$\La\subset\Fun(\Y)$, see Remark \ref{rmk1} below.

\begin{theorem}[cf. Theorem \ref{thmB}]
Let $z$, $z'$, and $\xi$ be complex parameters. For an arbitrary partition
$\nu=(\nu_1,\nu_2,\dots)$, there exists a unique element $\MM_\nu\in\La$, which
depends polynomially on $z,z'$ and rationally on $\xi$, and such that for any
natural number $N\ge\ell(\nu)$ and any $b>0$ and $\xi\in(0,1)$, one has
$$
\pi'_N\left(\MM_\nu\big|_{z=N,\,z'=N+b-1}\right)=M_{\nu\mid N,b,\xi}.
$$
\end{theorem}

\begin{definition}
We call the elements $\MM_\nu$ the {\it Meixner symmetric functions\/} with
parameters $z$, $z'$, and $\xi$.
\end{definition}

Next, we need the {\it Frobenius-Schur symmetric functions\/}. These are some
inhomogeneous elements $\FS_\nu\in\La$ indexed by arbitrary partitions $\nu$
and such that
$$
\FS_\nu=S_\nu\,+\,\textrm{lower degree terms}.
$$
For their definition, properties, and explicit expressions, see \cite{ORV03}.
In particular, one disposes of a simple explicit expression of the
Frobenius-Schur functions through the Schur functions.

\begin{theorem}[cf. Theorem \ref{thmC}]\label{thmC1}
The Meixner symmetric function $\MM_\nu$ with parameters $(z,z',\xi)$ are given
by the following expansion in the Frobenius-Schur symmetric functions:
$$
\MM_\nu=\sum_{\mu\subseteq\nu}C'(\nu,\mu;z,z',\xi)\FS_\mu,
$$
where
$$
C'(\nu,\mu;z,z',\xi)=(-1)^{|\nu|-|\mu|}
\left(\frac{\xi}{1-\xi}\right)^{|\nu|-|\mu|}\frac{\dim\nu/\mu}{(|\nu|-|\mu|)!}\,
\prod_{\Box\in\nu/\mu}(z+c(\Box))(z'+c(\Box))).
$$
\end{theorem}

\begin{definition}[cf. Theorem \ref{thmD}]
The {\it Meixner operator\/} $\D^\ME:\La\to\La$ with complex parameters
$(z,z',\xi)$ is defined in the basis $\{\FS_\nu\}$ of the Frobenius-Schur
functions by
$$
\D^\ME \FS_\nu=-|\nu|\FS_\nu+\frac\xi{1-\xi}
\sum_{\square\in\nu^-}(z+c(\square))(z'+c(\square))\FS_{\nu\setminus\square}.
$$
\end{definition}

{}From this definition one sees that $\D^\ME$ preserves the filtration of $\La$
and depends polynomially on $(z,z')$ and rationally on $\xi$, with the only
possible pole at $\xi=1$. The operator $\D^\ME$ is uniquely determined by these
properties together with the following one: If $z=N=1,2,\dots$, $z'=N+b-1$ with
$b>0$, and $\xi\in(0,1)$, then $\D^\ME$ preserves the kernel of the map
$\pi'_N:\La\to \La_N$ and the induced operator in $\La_N$ coincides with the
$N$-variate Meixner operator $D^\ME_N$ with parameter $b$.

\begin{theorem}[cf. Theorem \ref{thmE}]\label{thmE1}
The Meixner symmetric functions are eigenvectors of the Meixner operator $\D$
with the same values of parameters $(z,z',\xi)$. More precisely,
$$
\D^\ME\MM_\nu=-|\nu|\MM_\nu.
$$
\end{theorem}

\begin{definition}[cf. Definition \ref{dfnA}]\label{dfnC}
For any fixed $(z,z',\xi)\in\C^3$ with $\xi\ne1$, introduce the {\it formal
moment functional\/} $\psi^\ME:\La\to\C$ by setting
$$
\psi^\ME(1)=1, \qquad \psi(\D^\ME f)=0 \quad\textrm{for any $f\in\La$}.
$$
where the Meixner operator $\D^\ME$ is taken with the same values of the
parameters as $\psi$.
\end{definition}

\begin{theorem}[cf. Theorem \ref{thmF}]
For any Young diagram $\nu$,
$$
\psi^\ME(\FS_\nu)=\left(\frac{\xi}{1-\xi}\right)^{|\nu|}
\prod_{\Box\in\nu}(z+c(\Box))(z'+c(\Box))\cdot
\left(\frac{\dim\nu}{|\nu|!}\right)^2
$$
\end{theorem}

This formula provides an alternative (equivalent) way of introducing the moment
functional.

\begin{theorem}[cf. Theorem \ref{thmA}]\label{thmA1}
For any two Young diagrams $\nu$ and $\mu$,
$$
\psi^\ME\left(\MM_\nu\MM_\mu\right)
=\de_{\nu\mu}\frac{\xi^{|\nu|}}{(1-\xi)^{2|\nu|}}
\prod_{\Box\in\nu}(z+c(\Box))(z'+c(\Box)),
$$
where $\de_{\nu\mu}$ is Kronecker's delta.
\end{theorem}

This result shows that the Meixner symmetric functions are pairwise orthogonal
with respect to the inner product in the space $\La$ defined by
$$
(f,g)=(f,g)_{z,z',\xi}:=\psi^\ME(fg), \quad f,g\in\La.
$$
The inner product is nondegenerate provided that $z$ and $z'$ are not integers.

\begin{remark}
The two characterizations of the Laguerre symmetric functions from  Remark
\ref{rmk2} extend, with obvious modifications, to the Meixner symmetric
functions.
\end{remark}

\section{The orthogonality measure for the Meixner symmetric functions}

To speak about the orthogonality measure we have first to find an appropriate
realization of $\La$ as an algebra of functions on a space. In the context of
the Laguerre symmetric functions that space was the Thoma cone $\wt\Om$. Now
the relevant space is different: it is the countable set $\Y$ of Young
diagrams.

We will need the notion of {\it modified Frobenius coordinates\/} of a diagram
$\la\in\Y$. This is a double collection $(a;b)=(a_1,\dots,a_d;b_1,\dots,b_d)$
of half-integers, where $d$ stands for the number of diagonal boxes in $\la$,
$a_i=\la_i-i+\frac12$ equals the number of boxes in the $i$th row of $\la$ plus
one-half, and $b_i$ is the similar quantity for transposed diagram $\la'$. For
instance, if $\la=(3,2,2)$ then $(a;b)=(2\frac12, \frac12; 2\frac12,
1\frac12)$.

\begin{definition}
Let $\A$ be the unital algebra of functions on $\Y$ generated by the functions
of the form
$$
p_k(\la)=p_k(a;-b):=\sum_{i=1}^d[a_i^k -(-b_i)^k], \quad k=1,2,\dots,
$$
where $(a;b)$ is the collection of the modified Frobenius coordinates of a
diagram $\la\in\Y$. Elements of $\A$ are called {\it polynomial functions\/} on
$\Y$ \cite{KO94}.
\end{definition}

Consider the generators $p_1,p_2,\dots$ of $\La$ (the Newton power sums). The
assignment $p_k\mapsto p_k(\,\cdot\,)$ extends by multiplicativity to an
isomorphism $\La\to\A$ and hence defines an embedding of $\La$ into the algebra
$\Fun(\Y)$ of functions on the set $\Y$. This is the desired realization.

\begin{remark}\label{rmk1}
Now we can explain the origin of the map $\pi'_N:\La\to\La_N$ introduced in
Section \ref{s11}. Similarly to the realization $\La\subset\Fun(\Y)$, realize
$\La_N$ as an algebra of functions on $\Y_N\subset\Y$, the subset of Young
diagrams with at most $N$ nonzero rows, by letting the arguments $x_i$ of
$N$-variate symmetric polynomials to be equal to $\la_i+N-i$, where $\la$
ranges over $\Y_N$, $i=1\dots,N$. Then $\pi'_N$ is implemented by the natural
map $\Fun(\Y)\to\Fun(\Y_N)$ assigning to a function on $\Y$ its restriction to
$\Y_N$.
\end{remark}

For a diagram $\la\in\Y$, denote by $\la^+$ the set of the boxes that can be
appended to $\la$. As before, $\la^-$ is the set of the boxes that can be
removed from $\la$. By $\dim\la$ we denote the number of standard tableaux of
the shape $\la$.

\begin{theorem}
Under the realization $\La=\A\subset\Fun(\Y)$, the  Meixner  operator
$\D^\ME:\La\to\La$ with parameters $(z,z',\xi)$ is implemented by the following
operator in $\Fun(\Y)$, which will be denoted by the same symbol,
$$
\begin{aligned}
\D^\ME f(\la)&=\sum_{\square\in\la^+}A(\la,\square)f(\la\cup\square) +
\sum_{\square\in\la^-}B(\la,\square)f(\la\setminus\square)-C(\la)f(\la)\\
&=\sum_{\square\in\la^+}A(\la,\square)[f(\la\cup\square)-f(\la)]
+\sum_{\square\in\la^-}B(\la,\square)[f(\la\setminus\square)-f(\la)],
\end{aligned}
$$
where
\begin{equation*}
\begin{aligned}
A(\la,\square)&=\frac\xi{1-\xi}(z+c(\square))(z'+c(\square))
\frac{\dim(\la\cup\square)}{(|\la|+1)\dim\la}, \quad
\square\in\la^+,\\
B(\la,\square)&=\frac1{1-\xi}\sum_{\square\in\la^-}
\frac{|\la|\dim(\la\setminus\square)}{\dim\la}, \quad \square\in\la^-,\\
C(\la)&=\frac1{1-\xi}((1+\xi)|\la|+\xi zz').
\end{aligned}
\end{equation*}
\end{theorem}

\begin{definition}
Let $z$, $z'$, and $\xi$ be complex parameters,  $|\xi|<1$. The associated
complex measure on $\Y$, called the (mixed) {\it z-measure\/}, is defined by
\begin{equation*}
M_\zxi(\la)=(1-\xi)^{zz'}\prod_{\Box\in\La}(z+c(\Box))(z'+c(\Box))
\cdot\xi^{|\la|}\left(\frac{\dim\la}{|\la|!}\right)^2, \qquad \la\in\Y.
\end{equation*}
\end{definition}

One can prove that
$$
\sum_{\la\in\Y}M_\zxi(\la)=1.
$$
These measures were introduced in \cite{BO00}, and some closely related
measures on the finite sets of Young diagrams with a fixed number of boxes
appeared appeared earlier in \cite{KOV93}. The measures $M_\zxi$ are a special
case of Okounkov's Schur measures \cite{Ok01}.

The next theorem relates the measures $M_\zxi$ to the formal moment functional
$\psi^\ME$ with the same parameters (see Definition \ref{dfnC}).

\begin{theorem}[cf. Theorem \ref{thmG}]\label{thmG1}
Let $(z,z')$ be admissible and\/ $0<\xi<1$. Then the mixed z-measure $M_\zxi$
is a probability measure, all functions from $\La\subset\Fun(\Y)$ are
integrable with respect to $M_\zxi$, and
$$
\psi^\ME(f)=\sum_{\la\in\Y}f(\la)M_\zxi(\la), \qquad \forall
f\in\La\subset\Fun(\Y).
$$
Moreover, $\La$ is dense in the weight Hilbert space $\ell^2(\Y,M_\zxi)$.
\end{theorem}

This implies that (under the above assumptions on the parameters) $M_\zxi$
serves as the orthogonality measure for the Meixner symmetric functions.

\begin{remark}
The classical univariate Meixner polynomials are {\it autodual\/} in the sense
that, in an appropriate standardization, they are symmetric with respect to
transposition of the index and the argument, which both range over $\Z_+$. The
similar autoduality property holds for the Meixner symmetric functions viewed
as functions on $\Y$ under the realization $\La=\A\subset\Fun(\Y)$.

Indeed, under this realization, there is a simple expression for the functions
$\FS_\mu(\,\cdot\,)$ on $\Y$:
$$
\FS_\mu(\la)=\begin{cases}\dfrac{|\la|!}{(|\la|-|\mu|)!}\dfrac{\dim\la/\mu}{\dim\la},
& \textrm{if $\mu\subseteq\la$}\\
0, & \textrm{otherwise},
\end{cases} ,
$$
where $\la$ ranges over $\Y$; see \cite{ORV03}. Change the standardization of
$\MM_\nu$ by setting
$$
\MM_\nu=C''(\nu;z,z',\xi)\,\MM'_\nu,
$$
where
$$
C''(\nu;z,z',\xi):=\left(\frac{\xi}{1-\xi}\right)^{|\nu|}
\frac{\dim\nu}{|\nu|!} \prod_{\Box\in\nu}(z+c(\Box))(z'+c(\Box))
$$
is a normalizing factor. Then the above formula for $\FS_\mu(\la)$ combined
with Theorem \ref{thmC1} yields the following explicit expression for the
function $\MM'_\nu(\la)$:
\begin{multline*}
\MM'_\nu(\la)
=\sum_{\mu\subseteq(\nu\cap\la)}(-1)^{|\mu|}\left(\frac{1-\xi}{\xi}\right)^{|\mu|}
\frac{|\nu|!|\la|!}{(|\nu|-|\mu|)!(|\la|-|\mu|)!}\\
\times\frac{\dim\nu/\mu\,
\dim\la/\mu}{\dim\nu\,\dim\la}\prod_{\Box\in\mu}\frac1{(z+c(\Box))(z'+c(\Box))}.
\end{multline*}
Clearly, this expression is symmetric under $\nu\leftrightarrow\la$:
$$
\MM'_\nu(\la)=\MM'_\la(\nu), \qquad \nu,\,\la\in\Y.
$$
\end{remark}

\section{The Meixner jump process on the set of Young diagrams}

Return for a moment to the classical Meixner polynomials $M_n(x)$ and the
associated difference operator $D^\ME$ on $\Z_+$. Let $X_1^\ME(t)$ denote the
birth-death process $X_1^\ME(t)$ on $\Z_+$ whose jump rates are the
coefficients of $D^\ME$. That is, the rates of the jumps $x\to x+1$ and $x\to
x-1$ are equal to $(1-\xi)^{-1}\xi(b+x)$ and $(1-\xi)^{-1}x$, respectively.
This is a well-known instance of a birth-death process with linear jump rates.
The negative binomial distribution $w^\ME$ is the stationary distribution of
$X_1^\ME(t)$. The transition function of $X_1^\ME(t)$ can be expressed through
the Meixner polynomials according to formula \eqref{eq4}, where one has to
substitute $\phi_n=M_n$ and $w=w^\ME$.

More generally, the coefficients $A_i$ and $B_i$ of the operator $D^\ME_N$ (see
Section \ref{s11}) serve as the jump rates of a jump Markov process
$X^\ME_N(t)$ on the set $\Z^N_{+,\ord}$.

Even more generally, the following result holds (see \cite{BO06a}). Assume
$(z,z')$ is admissible and $0<\xi<1$. We know that then $M_\zxi$ is a
probability measure. Its support $\operatorname{supp} M_\zxi$ is the whole set
$\Y$ if $(z,z')$ belongs to the principal or complementary series, or a proper
subset of the form $\Y_N$ or $\{\la:\la'\in\Y_N\}$ if $(z,z')$ belongs to the
degenerate series.

\begin{theorem}
Let $(z,z')$ be admissible and $0<\xi<1$. Then there exists a jump Markov
process $X^\ME_\zxi(t)$ whose state space is $\operatorname{supp} M_\zxi$ and
whose jump rates are the coefficients $A(\la,\Box)$ and $B(\la,\Box)$ of the
Meixner operator $\D^\ME$. The measure $M_\zxi$ is an invariant and
symmetrizing measure for $X^\ME_\zxi(t)$.
\end{theorem}

The fact that the Meixner operator $\D^ME$ is diagonalized in the basis
$\MM_\nu$ of Meixner symmetric functions gives an expression for the transition
function, which we state for the case of nondegenerate parameters.

\begin{theorem}
Let $(z,z')$ belongs to the principal or complementary series and $0<\xi<1$.
The transition function $P(t;\la,\varkappa)$ of the Markov process
$X^\ME_\zxi(t)$ on $\Y$ can be written in the form form
$$
P(t;\la,\varkappa) =\sum_{\nu\in\Y}e^{-t|\nu|}
\frac{\MM_\nu(\la)\MM_\nu(\varkappa)}{\psi^\ME(\MM_\nu\MM_\nu)}\cdot
M_\zxi(\varkappa).
$$
Here $\la,\varkappa$ range over $\Y$, the $M_\nu$'s are viewed as functions on
$\Y$ in accordance with the realization $\La\subset\Fun(\Y)$, and the explicit
expression for $\psi^\ME(\MM_\nu\MM_\nu)$ is given in Theorem \ref{thmA1}.
\end{theorem}

\section{Approximation Meixner $\to$ Laguerre}\label{s14}

As well known, the Meixner polynomials $M_n(x)$ are discrete analogs of the
classical Laguerre polynomials $L_n(x)$. Namely, fix parameter $b>0$ and let
$\xi\uparrow1$. Then one has the limit relation
\begin{equation}\label{eq5}
\lim_{\xi\uparrow1}(1-\xi)^n M_n((1-\xi)^{-1}x)=L_n(x),
\end{equation}
where the scalar factor $(1-\xi)^n$ is used to keep the coefficient of $x^n$ to
be equal to 1.

The limit relation \eqref{eq5} can be easily derived from the explicit
expressions for the Laguerre and Meixner polynomials, see \cite{KS96}. On the
other hand, \eqref{eq5} can be explained by convergence of the weight
functions:  under the embedding
\begin{equation}\label{eq6}
\Z_+\to\R, \quad x\mapsto (1-\xi)x,
\end{equation}
the push-forward of the negative binomial distribution $w^\ME$ on $\Z_+$ with
parameters $(b,\xi)$ converges, as $\xi\uparrow1$, to the Gamma distribution
$(\Ga(b))^{-1}x^{b-1}e^{-x}dx$ on $\R_+$.

One more explanation can be given in terms of the univariate Meixner and
Laguerre operators. In the same scaling limit regime \eqref{eq6}, as the mesh
of the lattice goes to 0, the Meixner difference operator turns into the
Laguerre differential operator.

We are going to formulate similar statements in the infinite-dimensional
context.

Let $G:\La\to\La$ be the operator multiplying every homogeneous element by its
degree. In accordance with this, the operator $(1-\xi)^{-G}:\La\to\La$, which
appears in the next theorem, acts in the $m$th homogeneous component of $\La$
as multiplication by $(1-\xi)^{-m}$, for each $m\in\Z_+$.

The analog of \eqref{eq5} is

\begin{theorem}\label{thmH}
Let $\nu$ be an arbitrary partition, $\MM_\nu\in\La$ be the corresponding
Meixner symmetric function with parameters $(z,z',\xi)\in\C^3$, $\xi\ne1$, and
$\LL_\nu\in\La$ be the Laguerre symmetric function with parameters $(z,z')$.
One has
\begin{equation*}
\lim_{\xi\to1}(1-\xi)^{|\nu|} (1-\xi)^{-G}\MM_\nu=\LL_\nu,
\end{equation*}
where convergence holds in the finite-dimensional subspace of $\La$ consisting
of elements of degree $\le|\nu|$.
\end{theorem}

This is a direct corollary of Theorems \ref{thmC} and \ref{thmC1}.

As a corollary of Theorem \ref{thmH} combined with Theorems \ref{thmE} and
\ref{thmE1} one gets the following analog of the approximation of the
univariate Laguerre differential operator by the Meixner difference operator:

\begin{theorem}
Let\/ $\D^\ME:\La\to\La$ be the Meixner operator with parameters
$(z,z',\xi)\in\C^3$, $\xi\ne1$, and $\D:\La\to\La$ be the Laguerre operator
with parameters $(z,z')$. One has
$$
\lim_{\xi\to1}(1-\xi)^{-G}\circ\D^\ME\circ(1-\xi)^G=\D.
$$
\end{theorem}
Here we mean simple convergence on arbitrary elements $f\in\La$. Note that both
the pre-limit and limit operators preserve the filtration of $\La$, so that
application of the both operators to a given $f$ is contained in a fixed
finite-dimensional subspace of $\La$.

For $\epsi>0$ define the embedding
$$
\iota_\epsi:\Y\to\wt\Om, \quad \la=(a;b)\mapsto(\al,\be,\de),
$$
where $(a;b)=(a_1,\dots,a_d;b_1,\dots,b_d)$ are the modified Frobenius
coordinates of $\la\in\Y$, by setting
$$
\al_i=\begin{cases} \epsi a_i, & i\le d\\ 0, & i>d \end{cases}; \quad
\be_i=\begin{cases} \epsi b_i, & i\le d\\ 0, & i>d \end{cases}; \quad
\de=\epsi|\la|=\sum(\al_i+\be_i).
$$
The image $\iota_\epsi(\Y)$ is a discrete subset of $\wt\Om$. As
$\epsi\downarrow0$, it becomes more and dense in $\wt\Om$. This is the analog
of the embedding \eqref{eq6} ($\epsi=1-\xi$).

The analog of the approximation of the Gamma distribution by the negative
binomial distribution is

\begin{theorem}[cf. Section 5 in \cite{BO00}]
Let $(z,z')$ be admissible and $\xi\in(0,1)$. As $\xi\uparrow1$, the
pushforward of the measure $M_\zxi$ under the embedding
$\iota_{1-\xi}:\Y\to\wt\Om$ weakly converges to the measure $W_\zz$.
\end{theorem}

Recall that the weak topology on measures means convergence on bounded
continuous functions. Actually, more can be proved: convergence holds on any
test function on $\wt\Om$, which is continuous and grows, as
$\om=(\al,\be,\de)\to\infty$ not faster than a power of $\de$. In particular,
as test functions one can take elements of $\La\subset\Fun(\wt\om)$; then, by
virtue of Theorem \ref{thmG} and \ref{thmG1}, the claim of the theorem means
convergence of the formal moment functionals,
$$
\lim_{\xi\to1}\psi^\ME\circ(1-\xi)^G =\psi,
$$
which agrees with Theorems \ref{thmA} and \ref{thmA1}.

Finally, the approximation Meixner $\to$ Laguerre holds on the level of Markov
dynamics, which is used in the proof of the very existence of the diffusion
process $X(t)$.

\end{document}